\documentclass[]{amsart}
\def\O{{\mathcal O}}
\def\L{{\mathcal L}}
\def\P{{\mathbb P}}
\def\I{{\mathcal I}}

\def\p{{\mathbb P}}

\def\W{{K}}
\def\Pthree{{\mathbb P}^3}

\def\hx{{\widehat X}}
\def\wx{{\widetilde X}}
\def\wm{{\widetilde M}}
\def\wf{{\tilde f}}
\def\ws{{\tilde \sigma}}
\def\hp3{\widehat {\Pthree}}
\def\wp3{\widetilde {\Pthree}}
\def\ox{{\overline X}}
\def\om{{\overline M}}
\def\of{{\overline f}}
\def\os{{\overline \sigma}}
\def\op3{\overline {\Pthree}}

\def\Pone{{\mathbb P}^1}
\def\Aone{{\mathbb A}^1}

\def\Pic{\mathop{\rm Pic}}

\def\APic{\mathop{\rm APic}}

\def\Cl{\mathop{\rm Cl}}

\def\Spec{\mathop{\rm Spec}}
\def\codim{\mathop {\rm codim}}
\def\supp{\mathop {\rm Supp}}

\def\supp{\mathop{\rm Supp}}

\def\deg{\mathop{\rm deg}}

\def\cok{\mathop{\rm Coker}}
\def\Ker{\mathop{\rm Ker}}

\def\ss{\vskip .2 in}

\def\ra{\rightarrow}

\newtheorem{thm}{Theorem}[section]
\newtheorem{lem}[thm]{Lemma}
\newtheorem{claim}[thm]{Claim}

\newtheorem{prop}[thm]{Proposition}
\newtheorem{cor}[thm]{Corollary}
\newtheorem{defn}[thm]{Definition}

\newtheorem{rmk}[thm]{Remark}
\newtheorem{rmks}[thm]{Remarks}

\newtheorem{ex}[thm]{Example}
\newenvironment{pf}
   {\noindent {\bf Proof:}}{\vskip .20in}
\newenvironment{pfone}
   {\noindent {\bf Proof of Theorem \ref{one}: }}{\vskip .20in}

\date{}
\begin{document}

\title{Noether-Lefschetz Theorem with Base Locus}

\author{John Brevik}

\address{Dept. of Mathematics and Statistics, 
California State University,
Long Beach, CA 90840}

\email{jbrevik@csulb.edu}

\author{Scott Nollet}

\address{Mathematics Department, Texas Christian University, Fort 
Worth, TX 76129}

\email{s.nollet@tcu.edu}

\begin{abstract}
For an arbitrary curve $Z \subset \mathbb P^3$ (possibly reducible, non-reduced, unmixed) 
lying on a normal surface, the general surface $S$ of high degree containing $Z$ is also 
normal, but often singular. We compute the class groups of the very general such surface, 
thereby extending the Noether-Lefschetz theorem (the special case when $Z$ is empty). 
Our method is an adaptation of Griffiths and Harris' degeneration proof, simplified by a 
cohomology and base change argument. We give applications to computing Picard groups. 
\end{abstract}

\maketitle

\centerline{\it Dedicated to Robin Hartshorne on his 70th birthday}

\section{Introduction}

The algebraic surfaces in $\mathbb P^{3}_{\mathbb C}$ of degree $d$ are parametrized 
by the projective space $\mathbb P H^{0}(\O_{\Pthree} (d))$ via their equations. 
The Noether-Lefschetz locus $NL(d) \subset \mathbb P H^{0}(\O_{\Pthree} (d))$ 
corresponds to the smooth surfaces $S$ with $\Pic S \neq \langle \O_S (1) \rangle$. 
The original Noether-Lefschetz theorem, suggested by Noether in the 1880s and proved by 
Lefschetz in the 1920s, says that $NL(d)$ is a countable union of proper subvarieties if $d > 3$: 
in other words, the very general such surface $S$ satisfies 
$\Pic S = \langle \O_S (1) \rangle$. Work of Ciliberto, Green, Harris, Lopez, Miranda and Voisin 
around 1990 dramatically increased our understanding of the components of $NL(d)$. 
Mumford's challenge to write an explicit equation of a such a quartic surface $S$ was finally 
met by van Luijk \cite{vL} in the last few years. 

Carlson, Green, Griffiths and Harris proved the infinitesimal 
Noether theorem using variations of Hodge structures \cite{CGGH}; 
this method was used by Ein to compute the Picard group of the dependency locus 
of generic subspaces of sections of sufficiently ample vector bundles on 
arbitrary projective manifolds \cite{E} and by Green to explicitly bound the 
codimension of components of $NL(d)$ \cite{G2}. An approach from unpublished notes of 
Kumar and Srinivas led to Joshi's version for ambient projective threefolds 
\cite{J2} and the very recent extension to normal ambient varieties by 
Ravindra and Srinivas \cite{RS2}. 

While these recent generalizations are powerful and interesting, we have been 
impressed by Griffiths and Harris' degeneration method \cite{GH}, 
which relies on neither cohomological vanishings nor deformation theories, 
making it applicable to singular surfaces of low degree.
Lopez used it to compute Picard groups of general surfaces 
$S \subset \mathbb P^3$ containing a smooth connected curve $Z$ \cite[II.3.8]{L}, 
a result with many applications \cite{CL,CDE,FL}. We extend this result, 
replacing $Z$ with an arbitrary subscheme (possibly reducible, non-reduced, or 
of mixed dimension) which properly lies on a normal surface. Since the surfaces 
containing $Z$ are often necessarily singular, it is more natural to 
compute their class groups: 

\begin{thm}\label{main} 
Let $Z \subset \Pthree_{\mathbb C}$ be a closed subscheme of dimension 
$\leq 1$ with embedding dimension $\leq 2$ at all but finitely many
points and fix $d \geq 4$ with $\I_Z (d-1)$ generated by global sections.
Suppose that either
\begin{enumerate}
\item[(1)] $Z$ is reduced of embedding dimension two or  

\item[(2)] $H^0 (\I_Z (d-2)) \neq 0$. 
\end{enumerate}
Then the very general surface
$S \in |H^0 (\I_Z (d))|$ is normal with $\Cl S$ freely generated by 
$\O_{S} (1)$ and supports of the curve components of $Z$. 
\end{thm}

\begin{rmks}\label{hypoth}
    
(a) The hypotheses imply that $Z$ lies on a normal 
surface $T$ of degree $d-1$ with finitely generated class group. 
In case (1) this is because $Z$ actually lies on a smooth such surface 
and in case (2) this follows from Theorem \ref{one}. 
Thus Theorem \ref{main} follows from our more 
general Theorem \ref{genmain} and Proposition \ref{free}.
  
(b) With weaker hypothesis on $Z$, the general surface 
$S \in |H^0(\I_Z (d))|$ is not normal. Here one may consider the group
$\APic S$ of almost Cartier divisors from Hartshorne's theory of
generalized divisors \cite{GD}, but we would expect $\APic S$ 
to be infinitely generated due to the behavior along the curve part of the
singular locus \cite[Example 5.4 and Proposition 6.3]{GD}. 
\end{rmks}

As an application, we compute some Picard groups:

\begin{cor}\label{special} 
With the hypotheses of Theorem \ref{main}, let $Z_1, Z_2, \dots, Z_r$ 
be the curve components of $Z$. Then 
\begin{enumerate}
\item[(a)] If $\dim Z = 0$ (i.e. $r=0$), then $\Pic S = \langle \O_{S}(1) \rangle$.
\item[(b)] If $Z$ is a reduced l.c.i. curve and the $Z_i$ intersect at 
points of embedding dimension $2$ of $Z$, 
then $\Pic S = \langle \O_S (1), Z_1, \dots, Z_r \rangle$. 
\item[(c)] If $Z$ is an integral l.c.i curve, then $\Pic S = \langle \O_{S}(1),
\O_S (Z) \rangle$.
\item[(d)] If $Z$ has embedding dimension $2$, then 
$\Pic S = \langle \O_S (1), Z_1, \dots, Z_r \rangle$. 
\end{enumerate}
Moreover, the generating sets given in (b), (c), and (d) freely generate the given Picard groups. 
\end{cor}

\begin{pf} (a). Here $Z$ has {\it no} irreducible curve components, 
so $\Cl S = \langle \O_{S}(1) \rangle$, but $\O_{S} (1)$ is Cartier on $S$.
This strengthens Joshi's result extending the Noether-Lefschetz
theorem to very general {\it singular} surfaces \cite[4.4]{J}. 

(b) and (c). Here $\Cl S$ is generated as in the theorem, 
but the $Z_i$ intersect at smooth points of 
a general such surface $S$, so each $Z_i$ is Cartier on $S$ and we have
$\Cl S = \Pic S$.
The special case (c) where $r=1$ strengthens \cite[Cor. II.3.8]{L}.

For part (d), $\Cl S =  \langle \O_S(1), W_{1}, \dots W_{r} \rangle$ by
Theorem \ref{main}, where $W_i$ are the supports of the curve components
of $Z$. To compute $\Pic S$, we use the exact sequence
\begin{equation}\label{picloc}
0 \to \Pic S \to \Cl S \to \bigoplus_{\codim p = 2} \APic (\Spec \O_{S,p})
\end{equation}
introduced by Jaffe \cite{J} and developed by Hartshorne \cite[2.15]{GD}. 
The general surface $S$ is smooth where the components of $Z$ 
intersect and are singular at a finite number of points $p$ along the 
$Z_i$ of multiplicity $m_i > 1$,
each singularity having local equation $xy-z^{m_i}$ by 
Proposition \ref{singularities}(b).
For this type of singularity Hartshorne has shown that 
$\APic \Spec \O_{S,p} \cong \mathbb Z/{m_i \mathbb Z}$
generated by the class of $W_i$ \cite[Proposition 5.2]{Zeuthen}. 
Assembling the kernels of the pieces, we see 
that $\Pic S$ is generated by $\O (1)$ and $m_i W_i = Z_i$.
\end{pf}

\begin{ex}\label{fourlines}
The conclusion of Corollary \ref{special} (b) can fail if the $Z_i$ do not meet at points of embedding 
dimension $2$. The cone $Z \subset \Pthree$ over 4 planar points in general position consists of four 
lines $Z_i$ meeting at a point $p$ and is a complete intersection of 
two reducible quadrics, we may write $I_Z = (l_1 l_2, l_3 l_4)$. 
A very general degree $d \geq 4$ surface $S$ containing $Z$ is singular only at $p$ and has equation 
$F l_1 l_2 - G l_3 l_4 =0$ with 
$F(p), G(p) \neq 0$, hence the local ring of $S$ centered at $p$ can be written $\mathbb C[x,y,z]/(u l_1 l_2 - 
v l_3 l_4)$, where $u,v$ are units and 
$l_i$ are general linear forms. This is isomorphic to the local 
ring of the vertex of a quadric cone, so $\APic S \cong \mathbb Z/ 2 \mathbb Z$ 
generated by the class of a ruling \cite[II, Example 6.5.2]{AG} and any
of the lines $Z_i$ will do. Theorem. \ref{main} says that 
$\Cl S = \langle \O(1), Z_1,Z_2,Z_3,Z_4 \rangle$ and the map 
$\Cl S \to \APic S \cong \mathbb Z / 2 \mathbb Z$ in sequence (\ref{picloc})
takes $\O(1)$ to zero and each $Z_i$ to the generator, hence 
\[
\Pic S = \{\O(a)+\sum a_i L_i: 2|\sum a_i \} \subset \Cl S.
\]
\end{ex}

\begin{rmk}\label{franco}
One of our motivations is an application to Franco and Lascu's characterization of 
contractable curves \cite{FL}: if $Y \subset \mathbb P^3$ is an integral 
local complete intersection curve, then the following are equivalent:
\begin{enumerate}
\item[(1)] $Y$ is $\mathbb Q$-subcanonical.
\item[(2)] $Y$ is $\mathbb Q$-Gorenstein contractable to a point in a normal surface in $\mathbb P^3$ containing  
$Y$ as a Cartier divisor.  
\item[(3)] $Y$ is contractable on general surfaces of high degree. 
\end{enumerate}
The implication $(3) \Rightarrow (1)$ in their proof requires knowing that the Picard group of the general 
high degree surface $S$ containing $Y$ is generated by $Y$ and $\O_S (1)$, which is Corollary \ref{special} (b).
\end{rmk}

We also extend the Grothendieck-Lefschetz theorem for divisors with base
locus, using results of Ravindra and Srinivas \cite{RS}. First we note the weakest conditions that allow a 
subscheme to lie on a normal hypersurface.

\begin{defn}\label{superficial}
A closed subscheme $Z$ of a smooth ambient variety $M$ is 
{\em superficial} if (i) $\codim (Z,M) \geq 2$ and (ii) the closed subset 
$F \subset Z$ of points where the embedding dimension of $Z$ is equal to $\dim M$ satisfies $\codim (F,M) \geq 
3$.
\end{defn}  

\begin{thm}\label{one} Let $Z \subset \p^n_{\mathbb C}$ be superficial closed
subscheme with $n \geq 3$. If $\I_Z (d)$ is generated by global sections
and $H^0 (\I_Z (d-1)) \neq 0$, then 
\begin{enumerate}
\item[(a)] The Zariski general $S \in | H^0 (\I_Z (d)) |$ is normal with finitely generated class group. 
\item[(b)] If $n > 3$, then 
\begin{equation}\label{A}
\Cl S = \langle \O_S (1), W_1, \dots, W_r \rangle
\end{equation}
where $W_i$ are the supports of the codimension-$2$ components of $Z$. 
\end{enumerate}
\end{thm}

\begin{ex}\label{notfg}
The conclusions can fail if $H^0 (\I_Z (d-1)) = 0$: 
let $Z \subset \Pthree$ be the complete intersection of cones over two
smooth 
plane curves of degree $d>2$ with common vertex.  
The general surface $V$ of degree $d$ containing $Z$ is a cone over such a
curve $C$,
so $\Cl V \cong \Cl C = \Pic C$ \cite[II, Exer. 6.3 (a)]{AG} is infinitely
generated.
The cone over this example in $\p^4$ shows that part (b) also fails 
for $n>3$ without the $h^0 (\I_Z (d-1)) \neq 0$ condition.
\end{ex}

\begin{rmk}
A few words are in order comparing our main theorem to related results. 

When $Z$ is a smooth connected curve, theorem \ref{main} delivers 
exactly \cite[II, Corollary 3.8]{L} of Lopez.  Moreover, since the blow-up 
$M=\widetilde{\mathbb P^{3}} \to \mathbb P^{3}$ at $Z$ in diagram (\ref{std}) below 
embeds into $\mathbb P H^{0} (\I_{Z} (d))^{*}$ by strict transforms of the linear
system of degree-$d$ surfaces containing $Z$, and the general such surface maps
isomorphically onto its image in $\mathbb P^{3}$, this result also follows from 
results of Ein \cite{E}, Joshi \cite{J}, or the very recent theorem of Ravindra and 
Srinivas \cite{RS2} applied to $M$. 

The hypotheses for the results of Ein, Joshi, and Ravindra-Srinivas do not hold on $M$ 
for more general $Z$ (for example if $Z$ has some isolated points), so our result is of 
independent interest. 
On the other hand, these results apply to surfaces on more general
threefolds than $\mathbb P^3$, provided that some conditions (vanishing of cohomology,
respectively global generation of a certain sheaf) can be verified.
The main theorem of Lopez \cite[II. Theorem 3.1]{L} is also independent, 
as it relies on less restrictive hypotheses (the corresponding line bundle on $M$ 
need not be ample).
\end{rmk}

In section 2 we blow up the base locus of a linear system to interpret 
our results in terms of divisors on blow-ups and prove Theorem \ref{one}. 
Section 3 is an adaption of the degeneration method used by Griffiths, 
Harris and Lopez \cite{GH,L} to the case of families of singular surfaces.
For this we smooth the surfaces and form an \'etale cover of the family
where we can sort out the exceptional divisors. 
In the last section we prove the main theorem. Throughout we work over the
field $\mathbb C$ of complex numbers, since characteristic zero Bertini
theorems, generic smoothness, and monodromy arguments are used in the last two
sections. 
\ss
\noindent {\bf Acknowledgments:} As well as wishing him a happy birthday,
we thank Robin Hartshorne for his teachings, helpful comments and Example 
\ref{robin}. The second named author thanks Rosa Maria Mir\'o-Roig 
for asking a question which inspired this work and Andrew Sommese for
useful conversations.

\section{Finite generation of the class group}\label{standard}

Let $L$ be a line bundle on a smooth variety $M$ and 
$V \subset H^0 (M, L)$ a linear system defining a rational 
map $\phi: M \to \p V^*$. The image of the natural map 
$V \otimes L^{-1} \to \O_M$ defines the ideal of the base locus 
$Z \subset M$ for $V$. If $f:\wm \to M$ is the blow-up at $Z$, there is a closed
immersion 
$i: \wm \hookrightarrow M \times \p V^*$ whose image is the 
graph of $\phi$ and we have a diagram
\begin{equation}\label{std}
 \begin{array}{ccccc}
E & \subset & \wm & \stackrel{\sigma}{\to} & \p V^* \\
\downarrow & & f \downarrow & &\\
Z & \subset & M & & 
\end{array}
\end{equation}
with exceptional divisor $E$ and the map $\sigma = i \circ \pi_2$ given by the
invertible sheaf 
$\sigma^* (\O(1)) = f^*(L) \otimes \O_{\wm} (-E)$ 
(\cite[II, Example 7.17.3]{AG} and \cite[Theorem 1.3]{BS}). 

\begin{prop}\label{facts}
In the setting of diagram (\ref{std}), assume that $Z$ is superficial 
(Definition \ref{superficial}) with codimension-2 irreducible components $Z_i$.
Then 
\begin{enumerate}
\item[(a)] The general member $X \in |V|$ is normal. 
\item[(b)] Let $F \subset Z$ be the closed set where $\I_Z$ is not
2-generated or $Z$ has embedding dimension equal to $\dim M$. Then $\wm - f^{-1}(F)$ is normal
with class group
\[
\Cl (\wm - f^{-1}(F)) = \langle f^*(\Pic M), W_i \rangle
\]
where $W_i = {\overline {\supp f^{-1}(Z_i-F)}}$. 

\end{enumerate}
\end{prop}

\begin{pf} 
Since $Z$ is a codimension two local complete intersection off of $F$, 
the projection $\wm - f^{-1}(F) \to M-F$ is a $\Pone$-bundle 
over $Z - F$ and an isomorphism elsewhere. Letting $\Sigma \subset \wm$ 
denote the singular locus, $\Sigma - f^{-1}(F)$ is a set-theoretic section
over the non-smooth locus of $Z-F$ because the embedding dimension of $Z$ is less than $\dim M$ 
away
from $F$ \cite[Theorem 2.1]{N}, hence $\wm - f^{-1}(F)$ is regular in
codimension one (note that $Z$ superficial implies
 $\codim(F,M) > 2$ and $\codim(Z,M) \geq 2$).
Since $\wm - f^{-1}(F)$ is locally defined by a single equation in 
$(M - F) \times \Pone$, it satisfies Serre's $S_2$ condition and is
therefore normal. 
Moreover, each component 
$f^{-1}(Z_i)$ is supported on an irreducible Cartier divisor $W_i$
away from $\Sigma \cup f^{-1}(F)$ (because $\codim(F,M)>2$), hence 
\begin{equation}\label{ob2}
\Cl (\wm - f^{-1}(F)) = \langle f^*(\Pic M), W_i, \dots, W_r \rangle
\end{equation}
by repeated application of \cite[II, 6.5]{AG}.

For $X \in |V|$, view $\wx = f^{-1}(X)$ as a hyperplane section 
$\sigma^{-1}(H)$ with $H \in (\mathbb P V^*)^* = |V|$. 
Bertini theorems tell us that $\wx - f^{-1}(F)$ is regular in codimension
one. 
For $z \in Z_i-F$, the fibres $f^{-1}(z) \cong \Pone$ map isomorphically
to straight lines in $\p V^*$ because $-E$ is the relative $\O(1)$ in the
construction of the blow-up: since $f^{-1}(z)$ contains at most one
singular point of $\wm$ \cite[2.1]{N}, the general hyperplane 
$H$ meets this line transversely in a reduced point so the map 
$f:\wx \to X$ is a generic isomorphism along $Z$ (and an isomorphism away
from $f^{-1}(Z)$). Thus $X$ is regular in codimension one away from $F$,
and hence regular in codimension one because $\codim(F,M) > 2$. 
Divisors on smooth $M$ satisfy $S_2$, so general $X \in |V|$ are normal
\cite[II, Proposition 8.23 (b)]{AG}. 
\end{pf}

\begin{prop}\label{singularities} 
In the setting of Proposition \ref{facts} with $\dim M = 3$, the general surface 
$X \in |V|$ has singularities of two types: 
\begin{enumerate}
\item[(a)] {\em Fixed:} Points $F = \{z_j\}$ where $Z$ has embedding
dimension $3$.

\item[(b)] {\em Moving:} Away from $F$, there are a constant number of
singularities along each $Z_i$ with multiplicity $m_i > 1$; 
these move with $X$. For $X$ general, they have local equation 
$xy-z^{m_i}=0$.
\end{enumerate}
\end{prop}
\begin{pf}
Resuming the previous proof, we've seen that for $z \in Z_i - F$, the
general hyperplane $H \subset (\p V^*)^*$ meets the line 
$\sigma(f^{-1}(z))$ once, but we can say more. 
Consider the incidence  
\[
I=\{(z,H):z \in \bigcup Z_i - F, \sigma(f^{-1}(z)) \subset H\}.
\]
The fibres over the first projection to $\cup Z_i - F$ have dimension 
$\dim \p V -2$, so $\dim I = \dim \p V -1$ and 
$\dim \pi_2(I) \leq \dim \p V -1$, 
therefore the general hyperplane 
$H \in (\p V^*)^*$ meets {\it every} such line $\sigma(f^{-1}(z))$ once. 
It follows that $f:\wx - f^{-1}(F) \to X - F$ is an isomorphism in a
neighborhood of $\cup Z_i-F$, at least away from the isolated 
points of $Z$ which have embedding dimension two, where $X$ 
is already smooth. Therefore the singularities of $X$ away from $F$ 
are identified with those of $\wx$. 

If $G$ is the finite singular set of the support of $\cup Z_i - F$, 
$Z_i$ is locally a multiple of a smooth curve on a smooth surface 
away from $G$, so the ideal of $Z_i$ has the form $(x,y^{m_i})$, 
where $(x,y)$ is the ideal of the support, $m_i$ is the multiplicity, 
and $(x,y,z)$ is a regular 
sequence of parameters for the local ring $R=\O_{M,z}$. 
The blow-up of this ideal is covered by two affines, one being 
$\Spec R[u]/(ux-y^{m_i})$ which is singular exactly at the origin (the
other is smooth). 
We conclude that $\Sigma_i=\Sigma \cap f^{-1}(Z_i-F-G)$ is a section over
$Z_i-F-G$. 
If the image $\sigma(\Sigma_i)$ is a point, then $H$ misses 
$\sigma(\Sigma_i)$ and $\wx = \sigma^{-1}(H)$ is smooth along 
$f^{-1}(Z_i-F-G)$. 
If $\sigma(\Sigma_i)$ is an integral curve of degree $d_i$ and the map 
$\Sigma_i \to \sigma(\Sigma_i)$ has degree $e$, then $X^\prime$ has 
exactly $(e \cdot d_i)$ singularities along $\Sigma_i$ for general $H$. 
Since general $H$ meets $\sigma(\Sigma_i)$ transversely at each point, 
$\wx = \sigma^{-1}(H)$ has singularities with the same equation as above. 
\end{pf}

\begin{cor}\label{sing2} Let $Z \subset \mathbb P^3$ be superficial and
assume that $\I_Z (d)$ is generated by global sections. Then the general
surface $S$ of degree $d$ containing $Z$ is normal with constant number of
singularities, as described in Proposition \ref{singularities}.
\end{cor}

\begin{ex}\label{collapse}
For a concrete example, let $Z$ be 
the double structure on the line $L:x=y=0$ contained in the smooth cubic
surface $T \subset \Pthree$ with equation $x^{3}+y^{3}+xw^{2}+yz^{2}=0$ 
so that $I_{Z} = (x^{2},xy,y^{2},xw^{2}+yz^{2})$. A degree $d \geq 3$
surface $S$ 
containing $Z$ has equation 
\[
Ax^{2}+Bxy+Cy^{2}+H(xw^{2}+yz^{2})=0
\]
with $\deg A = \deg B =\deg C = d-2$ and $\deg H = d-3$: computing 
partial derivatives shows that this surface is singular at a 
point $q=(0,0,z_{0},w_{0})$ on $L$ precisely when $H(q)=0$, so for 
general $H$ there are exactly $d-3$ singular points. 

(a) If $d=3$, then the general surface $S$ is smooth since 
the special surface $T$ is. Here the constant number of moving 
singularities is $0$, even though $m_{1} > 1$. In this case 
$\sigma(\Sigma_{1})$ collapses to a point in the proof above. 

(b) For $d > 3$ and $H$ general, $S$ has exactly $d-3$ type-$A_{1}$ 
singularities along the line where $H =x=y=0$ as in 
Proposition \ref{singularities}.
For special $H$ meeting $L$ with higher multiplicity, these 
singularities can collide. 
\end{ex}

We close this section with a variant of the Grothendieck-Lefschetz
theorem for linear systems with base locus, which implies
Theorem \ref{one}. 

\begin{thm}\label{general} Let $L$ be a line bundle on a smooth projective 
variety $M$ and $V \subset H^0 (M, L)$ a linear system defining a rational
map $\phi: M \to \p V^*$ birational onto its image with superficial base
locus $Z$. Then the general $X \in |V|$ is normal and 
\begin{enumerate}
\item[(a)] If $\dim M > 3$, then 
$\cok(\Pic M \to \Cl X)$ is generated by the supports of the
codimension-$2$ components of $Z$.
\item[(b)] If $\dim M = 3$, then $\cok(\Pic M \to \Cl X)$ is finitely
generated.
\end{enumerate}
\end{thm}

\begin{pf}
Normality of $X$ is Proposition \ref{facts} (a). 
For the additional statements, let $\om \to \wm_{norm} \to \wm$ be the
normalization 
followed by a desingularization, with corresponding maps $\of, \wf$ 
to $M$ and $\os, \ws$ to $\p V^*$. Let $\ox$ (resp. $\wx_{norm}$) 
be a general hyperplane section of $\os$ (resp. $\ws$) and let 
$E_{\overline M}$ (resp. $E_{\overline X}$) be the union of exceptional
divisors for the 
desingularization $\om \to \wm_{norm}$ (resp. $\ox \to \wx_{norm}$). 
We have a commutative diagram
\[
\begin{array}{ccc}
\Pic \om & \stackrel{\rho}{\to} & \Pic \ox \\
\downarrow & & \downarrow \\
\Pic (\om - E_{\om}) & \to & \Pic (\ox - E_{\ox}) \\
\downarrow & & \downarrow \\
\Pic (\om - E_{\om} - {\of}^{-1}(F)) & \stackrel{\widetilde \rho}{\to} & 
\Pic (\ox - E_{\overline X} - {\of}^{-1}(F))
\end{array}
\]


Since $\tilde \sigma$ is birational, $\tilde \sigma^* \O(1)$ is a big
invertible sheaf, hence $\rho$ has finitely generated cokernel if $\dim M
\geq 3$ \cite[Theorem 2 (a)]{RS};
If $\dim M > 3$, the cokernel of $\rho$ is generated 
by divisors supported in $E_{\overline X}$ \cite[Theorem 2 (c)]{RS} so the 
middle horizontal map is surjective. Noting that the lower right 
vertical map is surjective, we conclude that $\widetilde \rho$ 
is surjective for $\dim M > 3$ and has finitely generated cokernel if 
$\dim M = 3$. 

Now because $\wm - f^{-1}(F)$ is normal (Proposition \ref{facts}(b)), the
desingularization  
$\om - {\of}^{-1}(F) \to \wm - f^{-1}(F)$ is obtained by blowing up smooth
centers in the singular loci (no normalization is required away from 
$f^{-1}(F)$) and we have the identifications 
$\om - E_{\om} - {\of}^{-1}(F)) \cong 
\wm - \Sigma - f^{-1}(F)$ 
and similarly $\ox - E_{\ox} - {\of}^{-1}(F) 
\cong \wx - \Sigma - f^{-1}(F)$.
Thus $\widetilde \rho$ may be identified with the restriction map $r$ 
\[
\Cl (\wm - f^{-1}(F))  \stackrel{r}{\to} \Cl (\wx - f^{-1}(F)).
\]
If $G \subset Z$ is the set over which $f: \wx-f^{-1}(F) \to X - F$ fails
to be an isomorphism (see proof of Proposition \ref{facts}), 
composing $r$ with the surjection
\[
\Cl (\wx - f^{-1}(F)) \to \Cl (\wx - f^{-1}(F) - f^{-1}(G)) 
\cong \Cl (X - F - G) \cong \Cl X
\]
shows that the last group is generated by $\Pic M$ and the supports of 
the $Z_i$ (use equation (\ref{ob2}) and note that the classes $W_i$ map to
$\supp Z_i$), so we draw conclusions (a) and (b). 
\end{pf}

\begin{pfone}
If $Z \subset \p^n$ is superficial, $\I_Z (d)$ is generated by global
sections 
and $H^0 (\I_Z (d-1)) \neq 0$, then $0 \neq f \in H^0 (\I_Z (d-1))$
implies that 
the rational map $\Pthree \to \p H^0 (\I_Z (d))^*$ is an isomorphism away
from 
the hypersurface $f=0$, hence birational onto its image and Theorem
\ref{general} applies. 
\end{pfone}

\section{Two Families of Surfaces}

We now turn to the harder problem of computing class groups of 
degree $d$ surfaces in $\Pthree$ with fixed base 
locus $Z$. In Proposition \ref{family1} we produce an open set 
$U \subset \p H^0 (\I_Z (d))$ and a family of desingularizations 
$\ox_t \to X_t$ for $t \in U$ such that the kernels of the natural 
maps $\Pic \ox_t \to \Cl X_t$ are represented by irreducible divisors 
over an \'etale cover $U^\prime \to U$. We restrict this family to a 
general pencil containing a reducible surface $T \cup P$, 
where (after modification) we compute
the Picard group of the central fibre (Proposition \ref{central}). 

Fix $Z \subset \Pthree$ superficial with $\I_Z (d-1)$ globally generated. 
Interpreting construction (\ref{std}) with $M= \Pthree$, $L =\O (d)$ 
and $V = H^0 (\I_Z (d)) \subset H^0 (L)$ yields a closed immersion 
$\sigma_{d}: \wp3 \hookrightarrow \p H^0 (\I_Z (d))^*$, where 
$\wp3 \stackrel{f}{\to} \Pthree$ is the blow-up at $Z$ \cite[Theorem
2.1]{BS}. Let $h: \op3 \to \wp3$ be a desingularization having smooth
exceptional divisors with normal crossings \cite{V} with composite maps 
$\os = \sigma_d \circ h: \op3 \to \p H^0 (\I_Z (d))^*$ and 
$\of = f \circ h: \op3 \to \Pthree$. 
There is a similar map $\sigma_{d-1}$ for degree $(d-1)$ surfaces which 
need not be a closed immersion.  

\begin{rmk}\label{planes}
If $P \in (\Pthree)^*$ is a general plane, we can describe its 
strict transform $\overline P \subset \op3$ and the map 
$\overline P \to P$. Following the proof of Proposition \ref{singularities}, 
let $F \subset Z$ be the finite set where $Z$ has embedding dimension 
three or $\I_Z$ is not 2-generated and let $G$ be the singularities of 
the support of $Z$. If $Z_i$ are the curve components of $Z$, then 
$Z_i$ has local ideal $(x, y^{m_i})$ away from $F \cup G$, where $m_i$ is
the multiplicity of $Z_i$. Therefore $\widetilde P \to P$ is the blow-up 
at $Z \cap P$ with exceptional divisors $W_{i,j}$ over the points in 
$Z_i$: these are the components of $W_i \cap \widetilde P$ with 
$W_i$ as in Proposition \ref{facts}. Moreover, $\widetilde P$ has exactly one
singularity in each $W_{i,j}$ with $m_i > 1$, which has equation $xy-z^{m_i}$. 
These singularities have canonical resolution compatible with the corresponding 
singular locus of $\wp3$ \cite[5.1 and 5.3]{Zeuthen} obtained 
by repeatedly blowing up points, so the fibres of $\overline P \to P$ 
over points in $Z_i \cap P$ consist of a connected chain of $m_i$ 
$\Pone$s, which include $W_{i,j}$. In particular, $\Pic \overline P$ is
freely generated by $\O(1), W_{i,j}$ and the exceptional divisors of the 
map $\overline P \to \widetilde P$ \cite[V, Cor. 5.4]{AG}.
\end{rmk}

We now compare the class groups of surfaces $S \in \p H^0 (\I_Z (d))$ 
and Picard groups of their strict transforms $\overline S \subset \op3$.
Letting $X \subset \Pthree \times \p H^0 (\I_Z (d))$ be the universal
family and $\ox \subset \op3 \times \p H^0 (\I_Z (d))$ be the family of
hyperplane divisors of $\os$, we have a diagram 
\begin{equation}\label{uni}
\begin{array}{ccc}
\ox & \subset & \op3 \times \p H^0 (\I_Z (d)) \\
\downarrow & & \downarrow \\
X & \subset & \Pthree \times \p H^0 (\I_Z (d)) \\
& & \downarrow \\
& & \p H^0 (\I_Z (d)).
\end{array}
\end{equation}

\begin{prop}\label{family1}
In the setting of diagram (\ref{uni}), there 
is a non-empty open set $U \subset \p H^0 (\I_Z (d))$, an 
\'etale cover $U^\prime \to U$, and effective irreducible divisors
$A_i \in \Pic \ox \times_U U^\prime$ 
such that 
\begin{enumerate}
\item[(a)] There is an open subset $V \subset U$ for which 
$\ox_V \to V$ is smooth and each surface $S_v$ with 
$v \in U^\prime \times_U V$ satisfies 
$\Ker (\Pic \overline S_v \to \Cl S_v) = \langle A_i \rangle$. 
\item[(b)] The set $U$ contains points corresponding to general 
reducible surfaces $X=T \cup P$ with $T \in \p H^0 (\I_Z (d-1))$ and 
$P \in (\Pthree)^*$. For these, $A_i \cap \overline T = \emptyset \iff A_i \cap \overline P \neq \emptyset$ 
and 
\begin{enumerate}
\item[1.] $\Ker (\Pic \overline T \to \Cl T) = \langle A_i \rangle$ 
\item[2.] $\Pic \overline P$ is freely generated by 
$\O(1)$, the strict transforms of the components $W_{i,j}$ of $W_i \cap \widetilde P \subset \wp3$ 
and the $A_i$ for which $A_i \cap \overline P \neq \emptyset$. 
\end{enumerate}

\end{enumerate}
\end{prop}

\begin{pf}
For $S \in \p H^0 (\I_Z (d))$, we view its strict transform 
$\overline S \subset \op3$ as a hyperplane section of the map $\os_d$, so 
$\overline S$ is smooth and irreducible by Bertini's theorem \cite{J1}. 
If $\Sigma \subset S$ is the singular locus, the kernel of the map 
\begin{equation}\label{picmap}
\Pic \overline S \to \Pic (\overline S - \of^{-1}(\Sigma)) \cong \Pic (S -
\Sigma) \cong \Cl S
\end{equation}
is generated by the irreducible divisors in $\of^{-1}(\Sigma)$, which appear as 
intersections with the following divisors in $\op3$: 
let $M_{i,j} \subset \op3$ be the irreducible divisors with 
$h(M_{i,j}) = \Sigma_i \subset \wp3$ 
(recall from Proposition \ref{singularities} that the singularities of $\wp3$ 
away from $f^{-1}(F)$ are sections $\Sigma_i$ of the curve components $Z_i 
\subset Z$, these give the moving singularities) and let $F_k \subset \op3$ be the 
irreducible divisors with $\of(F_k) \in F$ corresponding to fixed singularities. 
The kernels of the maps (\ref{picmap}) are generated by the components of 
$M_{i,j} \cap \overline S$ and $F_k \cap \overline S$. Similar statements 
apply to $T$ via the map $\os_{d-1}$.

Let $Q \subset \op3$ be a divisor $M_{i,j}$ or $F_k$ as described above. 

If $\dim h(Q)=0$, let $U_Q \subset \p H^0 (\I_Z (d))$ be the open subset of 
$H$ which miss $\sigma_d(h(Q))$. Clearly $U_Q$ is non-empty and contains 
reducible surfaces $T \cup P$, for general $H \in \p H^0 (\I_Z (d-1))$ misses 
$\sigma_{d-1}(h(Q))$ and general $P \in (\Pthree)^*$ misses $f(h(Q))$. 

If $\dim h(Q)=2$, let $U_Q \subset \p H^0 (\I_Z (d))$ be the non-empty open 
subset of $H$ for which $\os_d^{-1}(H) \cap Q$ is integral. Again $U_Q$ contains 
reducible surfaces $T \cup P$, for the general plane $P$ misses the point 
$f(h(Q)) \in F$. 
The closed immersion $\wp3 \hookrightarrow \Pthree \times \p H^0 (\I_Z (d-1))$
shows that $\sigma_{d-1}$ embeds $h(Q)$ into $\p H^0 (\I_Z (d-1))^*$, 
so Bertini's theorem tells us that $\overline T = \os_{d-1}^{-1}(H)$ 
meets $Q$ irreducibly. We will select $Q$ for one of the divisors $A_i$. 

Finally, if $\dim h(Q)=1$, then $h(Q)=C$ is an integral
curve and there is the Stein factorization of $Q \stackrel{h}{\to} C$
\[
Q \stackrel{\alpha}{\to} C^\prime \stackrel{\beta}{\to} C 
\stackrel{\sigma_{d}}{\hookrightarrow} \p H^0(\I_Z (d))^*
\]
in which the fibres of $\alpha$ are connected \cite[III, Cor. 11.5]{AG}. 
Since $Q$ is an exceptional divisor for $h$, it is smooth by construction 
and we may apply generic smoothness to the map 
$Q \stackrel{\beta \circ \alpha}{\to} C$ to find an open set 
$C^0 \subset C$ over which each fibre consists of exactly $d = \deg \beta$ smooth 
connected curves. Let $U_Q \subset \p H^0(\I_Z (d))$ be the open set of 
$H$ which meet $\sigma (C^0)$ in 
$e = \deg (C \stackrel{\sigma}{\to} \p H^0 (\I_Z (d)))$ reduced points, so 
that $\overline S \cap Q = \os^{-1}(H) \cap Q$ consists of exactly $de$ smooth 
connected curves. Then the lower horizontal map in diagram 
\begin{equation}\label{incidence}
\begin{array}{ccc}
{\mathcal Q} = \ox \cap (Q \times U_Q) & \to & U_Q \\
\downarrow & & \downarrow \\ 
I = \{(x^\prime,H): x^\prime \in C^\prime, H \in U_Q,
\sigma(\beta(x^\prime)) \in H\} & \to & U_Q \\
\downarrow & & \\
C^\prime & & 
\end{array}
\end{equation}
is an \'etale cover of degree $de$. 

We check that the open set $U_Q$ contains reducible surfaces $T \cup P$. 
If $Q = F_k$, then general $P$ misses $z$ and $\sigma_{d-1}$ embeds $h(Q)$ 
(as when $\dim h(Q)=2$ above). The map 
$\sigma_{d-1}: C \to \p H^0 (\I_Z (d-1))^*$ is given by the line bundle 
$f^*\O(d-1) \otimes \I_E$, where $E \subset \wp3$ is the exceptional divisor 
for the blow-up $f$, and $f^*\O(d) \otimes \I_E|_C \cong f^*\O(d-1) \otimes I_E|_C$ 
because $f^*(\O(d-1))|_C$ is trivial, so the map has degree $e$.
If $Q = M_{i,j}$, then $f(C) = W_i$ is the support of a curve component of $Z$ and 
the general plane $P$ meets $f(C)$ in $\deg W_i$ reduced points. In this case 
$C$ is a section of $W_i$, so $\deg f^* \O(1)|_C = \deg W_i$ and 
$\deg f^*\O(d) \otimes \I_E|_C = \deg f^*\O(d-1) \otimes \I_E|_C + \deg W_i$, 
so again $T \cup P$ will be in the degree $de$ \'etale locus $U_Q$. 

As an open subset of a projective bundle over $C^\prime$, $I$ is integral.  
To separate the $de$ connected components, base extend $I \to U$ by 
itself to obtain 
\[
\begin{array}{ccc}
I \times_U I & \to & I \\
\varphi \downarrow & & \downarrow \\
I & \to & U
\end{array}
\]
in which $\varphi$ is an \'etale cover of degree $de$ with the canonical diagonal 
section, so $I \times_U I$ is not connected. If $I^\prime \subset I \times_U I$ is 
any connected component for which the map 
$I^\prime \to I$ has degree $> 1$, we can base extend by 
$I^\prime \to I$ to split it up. 
We continue until we arrive at an integral base extension $U^\prime \to U$ 
for which the induced map $U^\prime \times_U I \to U^\prime$ is a trivial
\'etale cover of $de$ sheets. To finish, we base extend diagram
(\ref{incidence}) by 
$U^\prime \to U$. Since $I \times_U U^\prime$ has $de$ components,
so does 
${\mathcal Q} \times_U U^\prime \subset \ox \times_U U^\prime$, thereby
splitting the 
intersections $Q \cap \overline S$ into its components, which we will take 
as the $A_i$ in Proposition \ref{family1}. 

We carry out this procedure for each $Q$. Intersecting the resulting
Zariski open sets $U_Q$ and composing the finite \'etale covers, we 
obtain effective divisors $A_{i}$ on $\ox \times_U U^\prime$ which sort 
out the components of the intersections of the surfaces with the $Q$. 
Thus the $A_i$ generate the kernels of the maps 
$\Pic \overline S \to \Cl S$. To finish part (a), use generic smoothness 
to find an open subset $V \subset U$ where the map $\ox_V \to V$ is
smooth. Statement 2 of part (b) follows from Remark \ref{planes}. 
\end{pf}

\begin{rmk}\label{FG}
If $T_0 \in \p H^0 (\I_Z (d-1))$ is normal with finitely generated class
group, then a reducible surface $T \cup P$ as in Proposition \ref{family1}(b) 
can be chosen with $\Cl T$ finitely generated as well. 
If $T_1 \cup P_1$ is any surface in $U$, consider the 
linear deformation $T \to \Aone$ given by equation 
$(1-t) f_0 + t f_1=0$ in $\Pthree \times \Aone$, where $f_i=0$ is the
equation of $T_i$. Letting $\widetilde T \to T$ be a desingularization in 
which the central fibre $\widetilde T_0$ is smooth, the family 
$\widetilde T \to \Aone$ is flat and $\Pic \widetilde T_0$ is finitely
generated
because $\Cl T_0$ is, hence $H^1(\O_{\widetilde T_0})=0$. By
semicontinuity,
$H^1(\O_{\widetilde T_u})=0$ for $u$ near $0$, hence $\Cl T_u$ is finitely 
generated for $u$ near $0$. 
\end{rmk}

In computing the Picard group of a very general degree $d$ surface 
$S \subset \Pthree$ containing a smooth connected curve $Z$, 
Lopez \cite{L} adapted Griffiths and Harris' degeneration 
argument \cite{GH}. Our construction follows that of Lopez, 
except that we work in a blow-up of $\Pthree$ where the surfaces 
become smooth and must make a base extension to spread out divisors 
to compute the class groups. We extend
his \cite[Lem. II.3.3]{L} for these purposes.

\begin{lem}\label{inject}
Let $Z \subset \Pthree$ be superficial with curve components $Z_i$. 
Assume $\I_Z (d-1)$ is globally generated for some $d \geq 4$ and 
fix a normal surface $T \in \p H^0 (\I_Z (d-1))$ with finitely generated 
class group. Then the very general pair $(P,S) \in (\Pthree)^* \times \p
H^0 (\I_Z (d))$ 
with $D = T \cap P$ smooth satisfies
\begin{enumerate}
\item[(a)] The restriction map $\Cl T \to \Pic D$ is injective.
\item[(b)] $p \neq q \in Z_i \cap P \Rightarrow \O_D (p-q)$ is not torsion
in $\Pic D$ for each $i$.
\item[(c)] Let $p_{i,j}$ be the points in $Z_i \cap P$ and $q_k$ the 
remaining points in $(T \cap S) \cap P$. If $L \in \Cl T$ such that 
$L|_D \cong \O_D (\sum a_{i,j} p_{i,j} +\sum b_k q_k)$, 
then there are $\alpha_i$ and $\beta$ such that $a_{i,j}=\alpha_i$ for
each $j$ and $b_k = \beta$ for each $k$. 
\end{enumerate} 
\end{lem}

\begin{pf}
We follow the outline of \cite[Lem. II.3.3]{L}. 
First note that $T \subset \Pthree$ is not ruled by straight lines, for if
$T$ is 
a cone over a plane curve $C$, then normality of $T$ implies 
$C$ smooth and $\Cl T \cong \Pic C$ \cite[II, Example 6.3 (a)]{AG},
but the latter group 
is not finitely generated because $C$ is not rational. 
If $T$ is ruled but not a cone, then only finitely many rulings pass
through 
each singularity of $T$ and the general line $L \subset T$ is contained in 
the smooth locus. 
Here $\O_T (L)|_L \cong K_L \otimes K_T^\vee \cong \O_L (d-2)$ 
\cite[II,8.20]{AG}. 
The exact sequence $0 \to \O_T \to \O_T (L) \to \O_T (L)|_L \to 0$ shows
that 
$h^0 (L, \O_T (L)|_L) \leq 1$ (depending on whether $L$ is fixed on $T$ or
moves) 
and we conclude that $d < 3$, a contradiction.

For part (a), it is enough that $I(L) = \{D \in |\O_T(1)|: L_D = \O_D\}$ 
is a proper closed subset for $0 \neq L \in \Cl T$, since then the
countable union 
$\cup_{L \neq 0} I(L)$ cannot be all of $|\O_T(1)|$; thus, we show that
for fixed 
$L \in \Cl T$, $L|_D \cong \O_D$ for general $D \in |\O_T(1)|$ implies $L
\cong \O_T$. 

Because $T$ is not ruled by lines, the reducible plane sections in 
$|\O_T (1)|$ form a family of codimension $\geq 2$ \cite[Lemma II.2.4]{L}, 
so there is a pencil $\Pone \hookrightarrow {\p}H^0(\O_T(1))$ of irreducible 
curves whose base points lie in the smooth locus $T^0$. The total family 
\begin{equation}\label{total}
\begin{array}{ccccc}
{\widetilde T} & \subset & T \times \Pone & \stackrel{f}{\to} & \Pone \\
& & \downarrow g & & \\
& & T & & 
\end{array}
\end{equation} 
is isomorphic to the blow-up of $T$ at the base points and the exceptional
divisors 
map isomorphically onto $\Pone$ under $f$ by \cite[Theorem 1.3]{BS}.
Since ${\widetilde T} \stackrel{g}{\ra} T$ is an isomorphism near the
singularities, ${\widetilde L} = g^*(L)$ is reflexive on ${\widetilde T}$
and 
it suffices to show that ${\widetilde L} \cong \O_{\widetilde T}$.
The push-forward $B = f_*({\widetilde L})$ is a line bundle on $\Pone$,
for 
it is reflexive by \cite[Cor. 1.7]{SRS} and has rank $1$ because 
$h^0({\widetilde L}_t)=1$ for general $t \in \Pone$. 
Since $f_*({\widetilde L} \otimes f^*(B^\vee)) \cong \O_{\widetilde T}$ by 
the projection formula, there is 
$0 \neq s \in H^0({\widetilde L} \otimes f^*(B^\vee))$. 
The effective divisor $(s)_0$ does not map dominantly to $\Pone$ because 
${\widetilde L} \otimes f^*(B^\vee)|_t$ is trivial for general $t$. 
It follows that $(s)_0$ is a union of components of fibres of $f$, but 
since these are irreducible and $f_*(s)$ is nonvanishing, we conclude 
that $(s)_0 = \emptyset$, hence ${\widetilde L} \otimes f^*(B^\vee) \cong
\O_{\widetilde T}$ 
and $\widetilde L \cong f^*(B)$. 
If $x \in T$ is any base point for our pencil, 
$f^* B|_{g^{-1}(x)} \cong {\widetilde L}|_{g^{-1}(x)}$ is trivial, 
but $f:g^{-1}(x) \ra \Pone$ is an isomorphism, so $B= \O_{\Pone}$ and 
${\widetilde L} = \O_{\widetilde T}$.

For part (b), note that for $n > 1$ the set of planes $H$ for which there are
$p \neq q \in Z_{i} \cap H$
with $n \O_{D}(p-q)$ trivial is closed by semi-continuity, since this
condition is given 
by non-vanishing of a line bundle on a flat family; therefore it 
suffices to show the set of these planes is proper in $(\Pthree)^{*}$. 
For this we choose points $p,q \in Z_{i}$ such that the line $L$ through
$p,q$ meets $T$ at 
$d-1$ smooth points of $T$ and the general plane $H$ containing $L$ yields
a 
smooth curve $D = T \cap H$. As in part (a), the pencil of planes $H$ 
containing $L$ gives rise to a total family $\widetilde T$ isomorphic 
to the blow-up of $T$ at the points in $T \cap L$ and we again obtain 
diagram (\ref{total}). Let $E_{p}, E_{q} \cong \Pone$ be the 
exceptional divisors on $\widetilde T$ over $p,q$. 

The divisor $A_{n} = n \O_{\widetilde T} (E_{p}-E_{q})$ is non-trivial on 
$\widetilde T$ and restricts to $n \O_{D} (p-q)$ on the general fibre over 
$\Pone$: if the general such restriction is trivial on $D$, then the 
argument in part (a) shows that there is a line bundle $B \in \Pic \Pone$ 
such that $A_{n} \cong f^{*}(B)$. Since $d > 3$, there is a point 
$r \in L \cap T$ with $r \neq p,q$. The restriction of $A_{n}$ to 
$E_{r}$ is trivial, but $E_{r} \cong \Pone$ via the map $f$, so we 
see that $B$ itself is trivial and therefore $A_{n}$ is trivial on 
$\widetilde T$, a contradiction. 
We conclude that for each $n > 1$, $n \O_D (p - q)$ is only trivial for
finitely many $D$. 
Taking the union over $n > 1$ shows that the divisors $\O_D (p - q)$ are
not torsion for very general $H$. 

For part (c), we sketch Lopez' argument \cite[Lemma II.3.3 (3)]{L} 
with some improvements. Let $Y$ be the (integral) curve linked to $Z$ by 
$S \cap T$ and let
\[
I = \{(p_{i,j},q_k,P):\sum p_{i,j} + \sum q_k = P \cap (Z \cup Y)\}
\]
be the incidence set inside
\[
W_1^{d_1} \times W_2^{d_2} \times \dots \times W_r^{d_r} \times
Y^{d(d-1)-\deg Z} \times (\Pthree)^* 
\]
with projection $\pi$ onto the last component. 
Letting $U \subset (\Pthree)^*$ be the family of planes $P$ meeting $Z \cup Y$
transversely, 
$J = \pi^{-1}(U)$ is smooth and connected (the plane monodromy acts as the
product 
of symmetric groups \cite[Proposition II.2.6]{L}), so $J$ is irreducible.
For fixed $L \in \Cl T$ and $a_{i,j}, b_k \in \mathbb Z$, the sets
\[
J(L,a_{i,j},b_k)=\{(p_{i,j},q_k,P): \O_D(\sum a_{i,j} p_{i,j} + \sum b_k
q_k) \cong L|_D\} \subset J
\]
are closed by semicontinuity. Suppose that $J(L,a_{i,j},b_k)=J$.
Using the plane monodromy to permute the points $p_{i,j}$ for fixed $i$,
we arrive at the equation
$(a_{i,s} - a_{i,t})(p_{i,s}-p_{i,t})=0$ in $\Pic D$. By part (b), 
$J(L,a_{i,j},b_k)$ is a 
proper closed set if $a_{i,s} \neq a_{i,t}$ for any $s \neq t$, and so
has proper closed image in $(\Pthree)^*$. 
The countable union of all such images does not fill $(\Pthree)^*$, so for
very general $P$ we have 
$\O_D(\sum a_{i,j} p_{i,j} + \sum b_k q_k) \cong L|_D \Rightarrow a_{i,j}
= a_{i,j^\prime}, j \neq j^\prime$. 
For very general $S$ we have $b_k=b_{k^\prime}$ for $k \neq k^\prime$ 
because for $D$ fixed, we can vary $S$ to miss pairs $(p,q)$ with $\O_D
(p-q)$ torsion. 
\end{pf}

Now we construct the second family. Fix an integer $d \geq 4$, 
a superficial scheme $Z \subset \Pthree$ with $\I_Z (d-1)$ globally generated and assume that $Z$ lies 
on a normal surface of degree $d-1$. 
Letting $U \subset \p H^0 (\I_Z (d))$ be the open set constructed in 
Proposition \ref{family1}, Remark \ref{FG} tells us that there is a point $0 \in
U$ 
corresponding to a reducible surface $T \cup P$ with $T$ normal and 
$\Cl T$ finitely generated. Fix such a surface $T$ with equation $F=0$ 
and choose a plane $P \subset \Pthree$ with equation $L=0$ and a 
degree $d$ surface $S$ with equation $G=0$ as in Lemma \ref{inject} 
such that $T \cup P$ and $S$ are both in the set 
$U \subset \p H^0 (\I_Z (d))$ from Proposition \ref{family1}. 
These define the pencil of surfaces $S_t$ by equation $FL - tG =0$ for 
$t \in \Aone$. Set $\W = \{t \in \Aone: S_t \in U\}$. 
The universal property gives an embedding $\W \hookrightarrow U$.
Now base extend diagram (\ref{uni}) to obtain families

\begin{equation}\label{first}
\begin{array}{ccc} 
\ox_\W & \subset & \op3 \times \W \\
\downarrow & & \downarrow \\
X_\W & \subset & \Pthree \times \W \\
\end{array}
\end{equation}  

and an \'etale cover $\W^\prime \to \W$ with divisors $A_i$ as in 
Proposition \ref{family1} (a). 
If $\overline G, \overline F, \overline L$ are the local equations of
the strict transforms
$\overline S, \overline T, \overline P \subset \op3$, 
then the equation of $\ox$ is given by 
$\overline F \overline L - \overline G t= 0$, 
which is singular exactly at $t=0$ along the intersection 
$\overline S \cap \overline T \cap \overline P$. 
For $T,S$ general, Bertini assures us that $\overline Y = \overline S
\cap \overline T \subset \op3$ 
is a smooth connected curve, the strict transform of the curve 
$Y \subset \Pthree$ linked to $Z$ by the complete intersection $S \cap T$. 
Let $\hp3 \to \op3$ be the blow-up along $\overline Y$, giving an
associated family
\begin{equation}\label{third}
 \begin{array}{ccc} \hx & \subset & \hp3 \times \W \\
	& & \downarrow \\
	& & \W
    \end{array}
\end{equation}   
which agrees with $\ox_\W$ away from $t=0$ because $\overline Y$ 
is Cartier on $\ox_t$ for $t \neq 0$. 
At the central fibre $\widehat T \cong \overline T$ and 
$\widehat P \to \overline P$ is the blow-up along the reduced set of points 
$\overline P \cap \overline Y$. 
The resulting surfaces in $\hp3$ no longer intersect, so family
(\ref{third}) is smooth.

\begin{rmk}
Griffiths, Harris, and Lopez smoothed the family by blowing up the
quadratic singularities at the central fibre and blowing down the rulings
on the resulting quadrics. Local coordinate calculations show this is
equivalent to blowing up $\overline Y$ as above. 
\end{rmk}

We would like to compute $\Pic \hx_0$, but the monodromy of the moving
singularities causes ambiguity so instead we compute in the \'etale cover. 
Since $\W \subset U$, Proposition \ref{family1} gives an \'etale cover 
$e:\W^\prime \to \W$ and divisors $A_i$ on $\hx \times_\W \W^\prime$ 
which generate the kernels of the maps $\Pic \hx_t \to \Cl X_t$ for 
$t \neq 0$ and $\Pic \widehat T \to \Cl T$ at the central fibre. 
The exceptional divisor $\widehat Y$ for the blow-up 
$\hp3 \to \op3$ has the structure of a $\Pone$-bundle over 
$\overline Y$. Let $\widehat W_i \subset \hp3$ be the strict transforms of 
$W_i \subset \wp3$ from Observation \ref{facts} (b) 
and $\widehat P \subset \hp3$ be the strict transform of the plane $P$ 
supported in the central fibre. With this notation, we compute:

\begin{prop}\label{central}
In the setting of family (\ref{third}), let $e: \W^\prime \to \W$ be the
\'etale cover given in Proposition \ref{family1} and $A_i$ the corresponding
divisors on $\hx \times_\W \W^\prime$. 
For $p \in e^{-1}(0)$, we set 
$N = \O_{\hx \times_\W \W^\prime}(\widehat P)|_{(\hx \times_\W  
\W^\prime)_p}$, where $\widehat P$ is the irreducible component of 
$(\hx \times_\W \W^\prime)_p$ corresponding to $P$. Then
\begin{equation}\label{punchline}
\Pic (\hx \times_\W \W^\prime)_p = \langle \O (1), \widehat W_1, \dots,
\widehat W_r, A_i, N  \rangle.
\end{equation}
\end{prop}

\begin{pf}
The fibre $(\hx \times_\W \W^\prime)_p$ is the union 
$\widehat T \cup \widehat P$, where $\widehat T \cong \overline T$, 
$\widehat D = \widehat T \cap \widehat P$ is isomorphic to 
$\overline D \subset \overline T$ and $\widehat P \to \overline P$ 
is the blow-up at the $d(d-1)-\deg Z$ reduced points 
$\overline Y \cap \overline P$: let $Y_k \subset \widehat P$ be the
corresponding exceptional divisors. If $W_{i,j} \subset \widehat P$ are
strict transforms of the supports of the components of 
$W_i|_{\widetilde P}$ (Proposition \ref{facts}), then $\Pic \widehat P$ is 
freely generated by $\O(1), W_{i,j}, Y_k$ and the $A_i$ which meet 
$\widehat P$ by Remark \ref{planes} and Proposition \ref{family1} (b). Thus 
an arbitrary divisor 
\[
Q \in \Pic (\hx \times_\W \W^\prime)_p \cong 
\Pic \widehat T \times_{\Pic \widehat D} \Pic \widehat P 
\]
may be uniquely written as a pair
\[
Q=(A,\O(a) + \sum a_{i,j} W_{i,j} + \sum b_k Y_k + \sum_{A_i \cap \widehat P
\neq \emptyset} c_i A_i)
\]
with $A \in \Pic \widehat T$ and common restriction to $\widehat D$. 
Since $A_i|_{\widehat D}$ are trivial, tensoring
with $\O(-a)$ and applying Lemma \ref{inject} (c) shows that 
$a_{i,j}=\alpha_i$ and $b_k=\beta$ for some $\alpha_i$ and $\beta$
and $Q$ becomes  
\[
(A, \O(a) + \sum_{i,j} \alpha_i W_{i,j} + \beta \sum Y_k + \sum c_i A_i)=
(A, \O(a) + \sum \alpha_i \widehat W_i + \beta \widehat Y + \sum c_i A_i).
\]

Let us describe the divisor $N$ as a pair. Clearly $N_{\widehat T}$ 
is represented by the curve $\widehat D \in |\O_{\widehat T} (1)|$. 
The restriction of the divisor $\widehat T \cup \widehat P$ to itself is trivial, 
therefore $N_{\widehat P} = -{\widehat T}_{\widehat P}=-{\widehat D}$. 
As a divisor on $\widehat P$, $\widehat D$ takes the form 
$\O(d-1)-\sum m_i W_i - \sum Y_k - \sum_{A_i \cap \widehat P \neq
\emptyset} n_i A_i$ 
because the total transform of $D \subset P$ includes all the exceptional
divisors in the desingularization of $\widetilde P$. Since $Z_i$ meets $D$ at $m_i$-fold points, the 
supports $W_i$ have multiplicity $m_i$ in the total transform; the 
$A_i$ have positive multiplicity $n_i > 0$ whose exact values we will not need. 
Therefore 
\[
Q-\beta N = (A(-\beta), \O(a+\beta(1-d)) + \sum (\alpha_i +\beta m_i) \widehat W_i +
\sum_{A_i \cap \widehat P \neq \emptyset} (c_i-\beta) A_i).
\]
Finally, the kernel of $\Pic \widehat T \to \Cl T$ is generated by 
the $A_i$ meeting $\widehat T$ and $A(-\beta)|_{\widehat D}$ has form 
$\O(a+\beta(1-d)) + \sum (\alpha_i +\beta m_i) \widehat W_i$. 
Since $\Cl T \to \Pic D$ is injective by Lemma \ref{inject} (a), we 
conclude that 
\[
A(-\beta) = \O(a+\beta(1-d)) + \sum (\alpha_i +\beta m_i)
\widehat W_i + \sum_{A_i \cap \widehat T \neq \emptyset} d_i A_i
\]
and we have expressed $Q$ in terms of the generators stated. 
\end{pf}

\section{The Main Theorem}\label{degenerate}

Now we are now ready to prove Theorem \ref{main}. We define the 
relevant Noether-Lefschetz locus (\ref{bad}) and show that it is 
closed in the Hilbert-flag scheme (Proposition \ref{closed}).  
We show it is proper by reducing to families dominating the particular
family (\ref{third}) in Claim \ref{dominate}, where the proof is not
difficult. We begin with a useful consequence of cohomology and base
change.

\begin{prop}\label{windmill}
Let $S \stackrel{f}{\rightarrow} T$ be a projective flat family with $H^1 (\O_{S_t}) = 0$ for
each $t \in T$ 
and fix $L \in \Pic S$. Then 
\begin{enumerate}
\item[(a)] The set $\{t \in T: L_t = 0 \in \Pic S_t \}$ is open in $T$. 
\item[(b)] For $G \subset \Pic S$, the set 
$G_L = \{t \in T: L_t \in G_t \}$ is open in $T$.
\end{enumerate}
\end{prop}   

\begin{pf}
Suppose that $L_0 \cong \O_{S_0}$ for $0 \in T$. 
Then $H^1 L_0 = 0$ and this continues to hold in a 
Zariski open neighborhood $U \subset T$ about $0$ by semi-continuity.
Thus the natural map
$R^1f_* L \otimes k(t) \rightarrow H^1(S_t, L_t)=0$
is surjective over $U$, hence an isomorphism by 
cohomology and base change \cite[III, Theorem 12.11(a)]{AG}.
Therefore $R^1f_* L = 0$ is locally free on $U$, so again 
by cohomology and base change \cite[III, Theorem 12.11(b)]{AG} the
natural map 
\[  f_* L \rightarrow H^0(S_t, L_t) \]
is surjective and an isomorphism for all $t \in U$. 
Shrinking $U$ to an open affine if necessary, this gives 
surjectivity of $H^0(S,L)\rightarrow H^0(S_0,L_0)$,
allowing us to extend the nonvanishing global section  
$1 \in H^{0} \O_{S_{0}} \cong H^{0} L_0$ to a global 
section $s$ on $S_U$; since $s$ vanishes on a closed set, we can 
further shrink $U$ to avoid this set and obtain the desired result.
Part (b) follows because 
$G_L = \bigcup_{A \in G} \{t \in T: (L-A) = 0 \}$ is the union of open
sets. 
\end{pf}

\begin{ex}\label{robin}
Proposition \ref{windmill} (a) fails for $L \in \Cl S$ if we interpret 
$L_t = 0$ as $(L_t)^{\vee\vee} \cong \O_{S_t}$. 
For example, if $S_t$ is a family of smooth quadric surfaces degenerating
to the quadric cone $S_0 \subset \Pthree$ and 
$L \in \Cl S$ is the divisor which is the difference of opposite rulings
on $S_t$ for $t \neq 0$, then $L_t \neq 0$ for $t \neq 0$
but the limit is the difference of two rulings on the cone $S_0$, 
which is trivial. 
\end{ex}

Now we prove a stronger version of Theorem \ref{main}. We first prove 
that the proposed generators for the class groups have no relations. 

\begin{prop}\label{free}
Let $Z \subset \Pthree$ be superficial with curve components $Z_i$ having respective supports 
$W_i$ and assume $\I_Z (d-1)$ is generated by global sections for some $d \geq 4$.  
Fix open sets $V \subset U \subset \p H^0 (\I_Z (d))$ as in Proposition \ref{family1}. 
Then the divisors $\O(1), W_i$ have no relations in $\Cl X_t$ for $t \in V$. 
\end{prop}

\begin{pf}
Suppose that $\O(c) + \sum b_i W_i = 0$ in $\Cl X_t$ for some 
$t \in V$. Using the \'etale cover $V^\prime \to V$ from Proposition \ref{family1}, 
we lift $t$ to $t^\prime$ and have the same relation in $X_{t^\prime}$ and 
we obtain a corresponding relation $L=\O(c) + \sum b_i W_i + \sum a_i A_i=0$ in 
$\Pic \overline X_{t^\prime}$ (the $\sum a_i A_i$ terms are needed 
here because the pull-back of generalized divisors is not necessarily 
additive \cite[Example 2.18.1]{GD}).
Since the family is flat, $L_{t^\prime}$ is trivial in 
$\overline X_{t^\prime}$ on a Zariski open set $V_0 \subset V^\prime$ by 
Proposition \ref{windmill}. 

If $V_0$ is non-empty, then we can choose the surface $S$ from Proposition \ref{inject} to lie 
in $V_0$ so that $\W^\prime \cap V^\prime_R$ is non-empty. Thus the restriction of 
$L$ to $\overline X_{\W^\prime \cap V_0}$ is the trivial line bundle. 
By \cite[II, 6.5]{AG} it follows that $L$ extends to a line bundle on all of 
$\overline X_\W$, which is linearly equivalent to a combination of vertical components, 
hence we may write $\O(c) + \sum b_i W_i + \sum a_i A_i + \sum j V_j = 0$ on 
$\overline X \times_\W \W^\prime$ and after pulling back we have the same equation on 
$\widehat X \times_\W \W^\prime$. Restricting to $\widehat P$ in the central fibre, 
the vertical components away from the central fibre become trivial and the 
remaining combination of $\widehat P$ and $\widehat T$ 
is a multiple of $N$, so the relation in $\Pic \widehat P$ becomes  
$\O(c) + \sum b_i W_i + \sum a_i A_i + e N =0$. Now $\Pic \widehat P$ 
is freely generated by $\O(1), W_i, A_i$ and $Y_k$, so we see that $e=0$ because 
the coefficient of $Y_k$ in $N$ is nonzero (see calculation of 
$N|_{\widehat P}$ in the proof of Proposition \ref{central}), and therefore $c = b_i = 0$ as well.
\end{pf}

To finish the proof of Theorem \ref{main}, we must show that $\O(1)$ and the $W_i$ 
generate the class groups. First we show that the relevant locus in the Hilbert flag 
scheme is closed. Let
\[
{\mathcal H}_V = \{(C,S): C \subset S, S \in V \} 
\stackrel{\pi_2}{\to} V
\] 
be the Hilbert-flag scheme of locally Cohen-Macaulay curves $C$ on 
surfaces $S$ from the family $V$ and let
\begin{equation}\label{bad}
{\mathcal B}_V=\{(C,S) \in {\mathcal H}_V: \L(C) \not \in \langle \O_{S} (1),
W_{1}, 
\ldots W_{r} \rangle \subset \Cl S, S \in V \} 
\end{equation}
where $\L(C)=\I_{C,S}^{\vee}$ is the reflexive sheaf associated to $C$,
which generalizes the familiar line bundle $\O_S (C)$ on a smooth 
surface $S$ \cite{GD}. 

\begin{prop}\label{closed}
    ${\mathcal B}_V$ is closed in ${\mathcal H}_V$. 
\end{prop}

\begin{pf}
For each irreducible component ${\mathcal T} \subset {\mathcal H}_V$, we show that
${\mathcal B}_V \cap {\mathcal T}$ is closed in ${\mathcal T}$. Letting
\[
C  \subset  X  \subset  \Pthree_{\mathcal T}
\]
be the associated family of curves on surfaces, we base extend by a
desingularization to assume ${\mathcal T}$ is smooth and it suffices that 
\[
\{t \in {\mathcal T}: \L(C)_t \in \langle \O(1), W_1, \dots, W_r \rangle 
\subset \Cl X_t \}
\] 
is open in ${\mathcal T}$. 

If $V^\prime \to V$ is the \'etale cover from Proposition \ref{family1}, 
then ${\mathcal T} \times_V V^\prime \to {\mathcal T}$ is an \'etale cover of
smooth varieties. Since the total family $\ox_V \to V$ is smooth, 
so is $\ox \times_V ({\mathcal T} \times_V V^\prime)$ by base extension, 
and there are divisors $A_i$ which generate the kernels of the maps 
$\Pic {\overline X}_t \to \Cl X_t$. Now the pullback of $\L(C) \in \Cl X$ 
to $\ox \times_V ({\mathcal T} \times_V V^\prime)$ is a line bundle
$\overline \L$ and 
\[
\L(C)_t \in \langle \O(1), W_1, \dots, W_r \rangle\ \iff \overline \L_t \in
\langle \O(1), W_1, \dots, W_r, A_i \rangle.
\]
The latter condition is open by Proposition \ref{windmill}(b).
\end{pf}

\begin{thm}\label{genmain} 
Let $Z \subset \Pthree$ be a curve lying on a normal 
degree $(d-1)$ surface with finitely generated class group and assume that
$\I_Z (d-1)$ is generated by global sections. 
Then the very general surface $S \in |H^0 (\I_Z (d))|$ 
is normal with $\Cl S = \langle \O (1), W_1, \dots, W_r \rangle$. 
\end{thm}

\begin{pf}
Fix $V \subset \P H^0 (\I_{Z} (d)), {\mathcal H}_V$ and ${\mathcal B}_V$ as in
Proposition \ref{closed} above. The Hilbert-flag scheme ${\mathcal H}_V$ 
has countably many irreducible components, finitely many for each 
Hilbert polynomial for the curves $C$ in the family, hence 
${\mathcal B}_V$ also has finitely many irreducible components by 
Proposition \ref{closed}. To prove Theorem \ref{genmain} it suffices to 
show that if $D \subset {\mathcal B}_V$ is an irreducible component, then 
$\pi_2(D) \neq V$. Indeed, any surface $S \in V$ carrying a reflexive 
sheaf $L \not \in \langle \O(1), W_{1}, \dots, W_{r} \rangle$ also
contains such a curve $C$ (the zero section of $L(n)$ for some $n > 0$), 
so $(C,S) \in {\mathcal B}_V$ and $S \in \pi_2({\mathcal B}_V)$. 

To show that $\pi_2(D) \neq V$, it suffices to show that 
$\W \cap V \not \subset \pi_2(D)$, where $\W \subset \p H^0 (\I_Z (d))$ 
is the pencil of surfaces given by family (\ref{third}). 
If $\W \subset \overline {\pi_2(D)}$, we can apply Bertini's theorem 
to the closure $\overline {\pi_2^{-1}(\W)}$ in the full Hilbert scheme 
to obtain an integral {\it curve} $D$ dominating $\W$, so we may assume
$\dim D = 1$. Base extending by the normalization of $D$, it finally 
suffices to prove the following claim.

\begin{claim}\label{dominate} Let $X \subset \Pthree \times \W$ be the 
family (\ref{first}) and $f: D \to \W$ a surjective morphism of smooth
irreducible curves. 
If $C \subset Y = X \times_\W D$ is a flat family of smooth curves 
over $D$, then $C_{u} \in \langle \O(1), W_{1}, \dots, W_{r} \rangle$ 
for general $u \in D$.
\end{claim}    

Fix $p \in f^{-1}(0)$. Let $\widehat Y \to D$ be the base
extension of 
the family $\widehat X \to \W$ from (\ref{third}). We obtain
a diagram
\begin{equation}\label{family2}
    \begin{array}{ccccc}
    \widehat C & \subset & \widehat Y & \to & \widehat X \\
    & & \downarrow & & \downarrow \\
    & & D & \to & \W.
    \end{array}
\end{equation}
in which $\widehat C \subset \widehat Y$ is the strict transform of 
$C \subset Y$.
Letting $\widehat \L=\I_{\widehat C}^\vee \in \Cl \widehat Y$, we will
show that for $u$ near $p$,
$\widehat \L_u \in \langle \O(1), W_1, \dots W_r \rangle$.

{\underline {Case 1:}} If $f$ is unramified at $p$, then we have an
isomorphism 
$\widehat Y_p \cong \hx_0$ and the total family $\widehat Y$ is
smooth along 
$\widehat P \cap \widehat T = \widehat D$
in $\widehat Y_p$. After base extension by $\W^\prime \to \W$ we obtain
divisors $A_i$ as in Theorem \ref{family1}, and using Proposition \ref{central} 
(the hypotheses of these results holds because 
$Z$ lies on a normal surface, therefore is superficial) we may write 
\[
\widehat \L_p = \O (a) + b N + \sum a_i A_i + \sum b_j W_j
\]
Note that the hypotheses of Proposition \ref{windmill} apply to the family 
$\widehat Y \to D$ because the long exact cohomology sequence associated
to 
\[
0 \to \O_{\widehat Y_p} \to \O_{\widehat P} \oplus \O_{\widehat T} \to
\O_{\widehat D} \to 0
\]
shows that $H^1(\O_{\widehat Y_p})=0$, since 
$H^0 (\O_{\widehat D}) =\mathbb C$ and
$H^1 (\O_{\widehat T})=H^1(\O_{\widehat P})=0$. We apply it to the line 
bundle
\begin{equation}\label{trivial1}
{\mathcal M} = \O (a) + b N + \sum a_i A_i + \sum b_j W_j - \widehat \L \in
\Pic \widehat Y
\end{equation}
Since ${\mathcal M}_p$ is trivial, this is also true on an open neighborhood
of $p$ by Proposition \ref{windmill}(a).
The restrictions of $N = \O_{\widehat Y}(\widehat P)$ and $A_i$ are
trivial in $\Cl {\widehat Y}_u$ 
nearby, so $\widehat \L_u = \O (a) + \sum b_i W_i$ near $u=p$ and 
$C_u \in \langle \O(1), W_1, \dots, W_r \rangle$. 
\ss
{\underline {Case 2:}} Now suppose that $f$ is ramified at $p$.
We still have $\widehat Y_p \cong \hx_0$, but the ramification of
$f$ at $p$ causes the total family to be singular along $\widehat D$ in
the 
fibre $\widehat Y_p$. Specifically, let $y \in \widehat D \subset \widehat
Y_p$ 
be a point with image $z \in \widehat D \subset \hx_0$. Now $\hx$ is 
smooth at $z$ and if $(A,(t))=\O_{U,0}$ (resp. $(B,(u))=\O_{D,p}$) is the
local ring 
of $\W$ at $0$ (resp. $D$ at $p$), then the ring homomorphism 
$A \to B$ sends $t \mapsto u^s$ 
(up to unit) for some $s > 1$. Since $\hx$ is locally defined in
$\hp3 \times \W$ by the equation $L F - t = 0$, the base extension 
gives that $\widehat Y$ is locally 
defined in $\op3 \times D$ by $L F - u^s$. Thus locally speaking, 
$\O_{\widehat Y,y}$ is a quotient of a regular local ring $R$ in four
variables $L, F, H, u$ by the equation $L F - u^s = 0$. 

By \cite[5.1 and 5.3]{Zeuthen}, successive blowing up of the central curve
$\widehat D$ yields a desingularization $Z \to \widehat Y$ in which 
the fibre over $y$ is a chain of $\Pone$s and at the global level we have 
\[
Z_p = \widehat P \cup_{\widehat D_0} I_1 \cup_{\widehat D_1} \cup \dots
\cup I_{s-1} \cup_{\widehat D_{s-1}} \widehat T,
\]
where each $I_i$ is a ruled surface over both $\widehat D_i \cong
\widehat D_{i-1}$ and these two sections do not meet in $I_i$. 
Running exact sequences and induction shows that $H^1 (\O_{Z_p})=0$ as
before: 
for example, the exact sequence 
\[
0 \to \O_{\widehat P \cup I_1} \to \O_{\widehat P} \oplus \O_{I_1} \to
\O_{\widehat D_0} \to 0
\]
shows that $H^1 (\O_{\widehat P \cup I_1})=0$ because the induced map 
$H^1 (\O_{I_1}) \to H^1 (\O_{\widehat D})$ is an isomorphism via the 
section $\sigma:\widehat D \to I_1$. 
The total family $Z$ is smooth near the central fibre and (similar to
Proposition \ref{central}) we have
\[
\Pic Z_p = \langle \O (1), W_1, \dots W_r, A_i, N_0 = \O_Z (\widehat
P)|_{Z_p}, N_1, N_2, \dots N_{s-1} \rangle
\]
where $N_i = \O_Z (I_i)|_{Z_p}$ for $1 \leq i < s$, essentially because
every divisor on the ruled surface $I_i$ has the same restriction to 
$D_i$ and $D_{i-1}$ modulo $\O_{I_i} (D_i)$.

Now the proof goes through as in the unramified case, the point being that
the new divisors $N_i$ have trivial restrictions in the nearby $\Cl Z_u$. 
\end{pf}

\end{document}